\documentclass[a4paper, 10pt]{article}
\usepackage{amsmath,amsthm,amssymb}
\usepackage[T2A]{fontenc}
\usepackage[utf8]{inputenc}
\usepackage[english, russian]{babel}
\usepackage{graphicx}
\usepackage{enumerate}
\usepackage{tikz}
\usepackage{geometry}
\begin{document}
\begin{center}
\LARGE\textbf{Multiplier Theorem for Fourier Series in continuous-discrete Sobolev orthogonal polynomials.}

\large \textbf{(B.P.Osilenker, Moscow, b\_osilenker@mail.ru)}
\end{center}

\setlength{\leftskip}{5em}
\setlength{\rightskip}{5em}
\textbf{Abstract.} In paper we study the multipliers of Fourier series in polynomials orthogonal in continuous-discrete  Sobolev’s spaces. Multiplier Theorem for Fourier-Sobolev series is obtained.This result  based on the representation of the Fejér kernel, on the construction of  the “humpbacked majorant” and weighted estimates of maximal functions.


\setlength{\leftskip}{0em}
\setlength{\rightskip}{0em}

Let $\theta(x)$ be a positive Borel measure in $[-1, 1]$, with infinitely many points of increasing and let the masspoints $a_k$, $-1\leq a_k \leq 1$, $k=1,2,\ldots, m$. For $f$ and $g$ in $L^2_{d\theta}[-1, 1]$ such that there exist the derivatives in $a_k$, one can introduce the inner product

\begin{equation}
\left<f, g\right> = \int^1_{-1} f(x)g(x)d\theta(x) + \sum^m_{k=1}\sum^{N_k}_{i=0} M_{k, i} f^{(i)}(a_k)g^{(i)}(a_k), \nonumber
\end{equation}
where $M_{k, i}\geq 0\ (i=0,1,2,\ldots,N_k-1);\ M_{k, N_k}>0, k=1,2,\ldots,m)$ and $\theta(\{a_k\})=0\ (k=1,2,\ldots,m)$.

Linear spaces with this inner product is called «continuous-discrete Sobolev spaces».

Let $\{\widehat{q}_n(x),\ n\in\mathbb{Z}_+,\ \mathbb{Z}_+=\{0,1,2,\ldots\};\ x\in[-1, 1]\}$ be the sequence of polynomials of degree $n$ with a positive leading coefficients orthonormal with respect to the inner product(continuous-discrete Sobolev orthonormal polynomials)

\begin{equation}
\left<\widehat{q}_n, \widehat{q}_m\right> = \int^1_{-1} \widehat{q}_n(x)\widehat{q}_m(x)d\theta(x) + \sum^m_{k=1}\sum^{N_k}_{i=0} M_{k, i} \widehat{q}_n^{(i)}(a_k)\widehat{q}_m^{(i)}(a_k) = \delta_{n,m}. \nonumber
\end{equation}

Denote by $\mathfrak{R}_p(1\leq p <\infty)$ the set of functions
\begin{equation}
\mathfrak{R}_p = \begin{Bmatrix}
f,\ \int_{-1}^1|f(x)|^pd\theta(x)<\infty;\ f^{(i)}(a_k)\text{ - exist}\\
i = 0,1,2,\ldots,N_k;\ -1\leq a_k\leq 1(k=1,2,\ldots,m)
\end{Bmatrix}. \nonumber
\end{equation}

To each $f\in\mathfrak{R}_p$ we assign Fourier-Sobolev series
\begin{equation}
f(x) \sim \sum^\infty_{k=0}c_k(f)\widehat{q}_k(x)\ (x\in[-1, 1]), \nonumber
\end{equation}

with Fourier coefficients
\begin{equation}
c_k(f) = \left<f, \widehat{q}_k\right> = \int^1_{-1}f(x)\widehat{q}_k(x)d\theta(x) +\sum^m_{s=1}\sum^{N_s}_{i=0} M_{s,i}f^{(i)}(a_s)\widehat{q}^{(i)}_k(a_s). \nonumber
\end{equation}

We consider the following sequence of the real numbers
\begin{equation}
\Phi = \{\phi_k,\ k=0,1,2,\ldots;\ \phi_0=1;\ \{\phi_k\}^\infty_{k=0}\in 1^\infty\}. \nonumber
\end{equation}

For any function $f\in\mathfrak{R}_p$ by their Fourier-Sobolev series we introduce the linear transformation $T$ defined by relation
\begin{equation}
\label{eq1}
T(f;x;\Phi)\sim\sum^\infty_{k=0}\phi_k c_k(f)\widehat{q}_k(x).
\end{equation}

Transformation $T$ is called the multiplier operator, the sequence $\{\phi_k\}^\infty_{k=0}$ is called the multiplier of convergence6 and series \eqref{eq1} is called the multiplier of convergence6 and series \eqref{eq1} is called the multiplier series.

We investigate some problems of pointwise and uniform multipliers of convergence for Fourier-Sobolev series. Multiplier Theorem for the Fourier-Sobolev series is obtained.
There are many papers have been devoted to continuous-discrete Sobolev orthonormal polynomials and Fourier series (see, for example [1] –[27]).

Some results about multipliers of the Fourier series in polynomials orthonormal in continuous-discrete Sobolev spaces were announced in [20].

Let $N^*_k$ be-the-positive integer number defined by
\begin{equation}
N^*_k = \begin{cases}
N_k + 1,\text{ if }N_k\text{ is odd,}\\
N_k + 2,\text{ if }N_k\text{ is even,}
\end{cases}\nonumber
\end{equation}
\begin{equation}
w_N(x) = \Pi^m_{k=1} (x-a_k)^{N^*_k},\ N=\sum^m_{k=1}N^*_k,\ \pi_{N+1}(x)=\int^x_{-1} w_N(t)dt. \nonumber
\end{equation}

Orthonormal polynomials $\widehat{q}_n(x)$ satisfy the following recurrence relation
\begin{equation}
\pi_{N+1}\widehat{q}_n(x) = \sum^{N+1}_{j=0} d_{n+j,j}\widehat{q}_{n+j}(x) + \sum^{N+1}_{j=1} d_{n,j}\widehat{q}_{n-j}(x)(n\in\mathbb{Z}_+;\widehat{q}_{-j}=0,\ j=1,2,\ldots;\ d_{n,s}=0,\ n=0,1,\ldots,s-1).   \nonumber
\end{equation}

Define by
\begin{equation}
\varepsilon_m = (-1, 1)\ \cup^m_{s=1}\{a_s\}. \nonumber
\end{equation}

The sequence $\Phi = \{\phi_n,\ n=0,1,2,\ldots;\ \phi_0=1\}$ is called quasiconvex if
\begin{equation}
\sum^\infty_{k=0} (k+1)|\Delta^2\phi_k|<\infty,\nonumber
\end{equation}
where $\Delta\phi_k = \phi_k-\phi_{k+1}$, $\Delta^2\phi_k = \Delta(\Delta\phi) = \phi_k - 2\phi_{k+1} + \phi_{k+2}(k=0,1,\ldots,n)$.

\textbf{Theorem 1.} Let the orthonormal polynomial system $\{\widehat{q}_k(x)\}^\infty_{k=0}$ be satisfy the following condition
\begin{equation}\label{eq2}
|\widehat{q}_k(t)|\leq h(t)(t\in\varepsilon_m)
\end{equation}
and for the recurrence coefficients the estimate
\begin{equation}\label{eq3}
\sum^{N+1}_{j=1}j\sum^{N+1}_{l=0}\sum^\infty_{s=0}\left(|d_{s+j,j}-d_{s+j+l,j}| + |d_{s+j,l}-d_{s+j+l,l}|\right)<\infty
\end{equation}
holds. If for quasiconvex sequence $\Phi$ the relation
\begin{equation}\label{eq4}
\phi_k = O\left(\frac{1}{\ln k}\right)(k\rightarrow\infty)
\end{equation}
holds, then the following statements are valid:
\begin{enumerate}[(i)]
\item let for each function $f\in\mathfrak{R}_p(1\leq p<\infty)$ be fulfilled
\begin{equation}\label{eq5}
\int^1_{-1}\left|f(t)\right|^ph^p(t)d\theta(t)<\infty,\ \ \int^1_{-1}h^p(t)d\theta(t)<\infty,
\end{equation}
then at every Lebesgue's point $x\in\varepsilon_m$(and, consequently, a.e.) the series \eqref{eq1} converges
\begin{equation}
T(f;x;\Phi) = \sum^\infty_{k=0}\phi_kc_k(f)\widehat{q}_k(x);\nonumber
\end{equation}
\item in addition, suppose function $f$ is continuous in $[-1,1]$ and the measure $d\theta(x)$ is absolutely continuous and
    \begin{equation}\label{eq6}
    d\theta(x) = \omega(x)dx,\ \omega(x)\text{ is continuous in }\varepsilon_m;
    \end{equation}
    then the series \eqref{eq1} is uniformly converges on compact subsets $K\subset\varepsilon_m$.
\end{enumerate}

We define for $f\in\mathfrak{R}_p$ the space $W^p_\theta(F)(1\leq p<\infty)$ for subset $F\subseteq[-1, 1]$:
\begin{equation}
W^p_\theta(F) = \{f,\ ||f||_{W^p_{\theta}(F)}<+\infty,\ ||f||^P_{W^p_{\theta^p}(F)} = ||f||^p_{L^P_\theta(F)} +\sum^m_{k=1}\sum^{N_k}_{i=0}M_{k,i}|f^{(i)}(a_k)|^p\}. \nonumber
\end{equation}

The space $W^p_\theta([-1,1])(1\leq p <\infty)$ is not complete.

\textbf{Theorem 2.} Let the orthonormal polynomial system $\{\widehat{q}_k(x)\}^\infty_{k=0}$ be satisfy the following condition \eqref{eq2}, \eqref{eq3}, \eqref{eq6} and
\begin{equation}\label{eq7}
\sup_{n\in\mathbb{Z}_+}\sum^n_{j=0}|q^{(i)}_j(a_s)|<\infty(i=0,1,\ldots,N_s;\ s=1,2,\ldots,m),
\end{equation}
\begin{equation}\label{eq8}
||h||_{L^p_\theta([-1,1])}<\infty,\ ||h||_{L^q_\theta([-1,1])}<\infty\left(1<p<\infty,\ \frac{1}{p}+\frac{1}{q}=1\right).
\end{equation}

If the sequence $\Phi$ is quasiconvex and satisfy \eqref{eq4}, then for $f\in W^p_\theta([-1,1])(1<p<\infty)$, satisfying \eqref{eq5}, on any compact subsets $K\subseteq \varepsilon_m$ the following estimate
\begin{equation}
||T(f;x;\Phi)||_{W^p_\theta(K)}\leq C_p||f||_{W^p_\theta([-1,1])},\nonumber
\end{equation}
holds, where the constant $C_p>0$ indepent on function $f$ and the sequence $\Phi$.

\textbf{Remark.} Symmetric Gegenbauer-Sobolev orthonormal polynomials $\left\{\widehat{B}^{(\alpha)}_n(x)\right\}$$\left(n\in\mathbb{Z}_+;\ x\in[-1,1]\right)$ orthonormal in an inner product
\begin{equation}
\left<f,g\right>_\alpha = \int^1_{-1} f(x)g(x)w_\alpha(x)dx + M\left[f(1)g(1) + f(-1)g(-1)\right] + N\left[f'(1)g'(1) + f'(-1)g'(-1)\right]\ \ \ (M\geq0;N\geq0),\nonumber
\end{equation}
where
\begin{equation}
w_\alpha(x) = \frac{\Gamma(2\alpha+2)}{2^{2\alpha + 1}\Gamma^2(\alpha+1)}(1-x^2)^\alpha\left(\alpha>\frac{1}{2}\right),\nonumber
\end{equation}
satisfying the conditions \eqref{eq2}, \eqref{eq3}, \eqref{eq6}, \eqref{eq7}, \eqref{eq8}.

\begin{center}
\large \textbf{REFERENCES}
\end{center}
\begin{enumerate}[{[}1{]}]
\item H.Bavinck. Differential operators having Sobolev-type Gegenbauer polynomials as eigenfunctions, J.Comput. Appl. Math., 118(2000), 23-42.
\item H.Bavinck, J.Koekoek. Differential operators having symmetric orthogonal polynomials as eigenfunctions,  J.Comput. Appl. Math.,106(1999), 369-393.
\item H.Bavinck, H.G.Meijer.Orthogonal polynomials with respect to a symmetric inner poduct involving   derivatives, Appl. Anal. 33(1989), 103-117.
\item H.Bavinck, H.G.Meijer .On orthogonal polynomials with respect to an inner product involving derivatives: zeros and recurrence relations, Indag. Math.(N.S.) 1(1990), 7-14.
\item Ó.Ciaurri, J.Minguez Fourier series of Gegenbauer-Sobolev polynomials. Symmetry,Integrability, geometry, Methods and Applications 14(2018), Paper 024, 11pp.
\item Ó.Ciaurri, J.Minguez Fourier series of   Jacobi-Sobolev polynomials. Transforms Spec. Funct., 30(2019), 334-346.
\item Foulquie A.Moreno,F. Marcellán, B.P. Osilenker, Estimates for polynomials orthogonal  with respect to some Gegenbauer-Sobolev inner product, J.Ineq.Appl. 3(1999), 401-419.
\item Abel Díaz-González, Francisco Marcellán-Español , Héctor Pijeira-Cabrera and Wilfredo Urbina-Romero. Discrete-Conti\-nuous Jacobi-Sobolev Spaces and Fourier Series. arXiv:1911, 12746v1[math.CA]28 Nov.2019.
\item J. Heinonen, T.Kilpelâinen, O.Martio, Nonlinear Potential Theory of Degenerate Elliptic Equations, Oxford Science Publ., Clarendon Press, Oxford, 1993.
\item T.Kilpelâinen, Weighted Sobolev spaces and capacity, Ann. Acad. Sci. Fenn. Ser.A.I.Math. 19(1994), 95-113.
\item R.Koekoek. Differential equations for symmetric generalized ultraspherical
polynomials//Trans.Amer.Math.Soc. 345:1 (1994), 47-72.
\item A.Kufner, Weighted Sobolev Spaces, Teubner Verlagsgesellschaft, Teubner – Texte  zur Mathematik (Band 31), 1980; also published by Wiley, 1985.
\item A.Kufner,A.M.Sänding, Some Applicatios of Weighted Sobolev Spaces, Teubner Verlagsgesellschaft, Teubner – Texte  zur Mathematik (Band 100), [39].
\item F.Marcellán, B.P.Osilenker, I.A.Rocha, On Fourier series of Jacobi- Sobolev orthogonal polynomials, J. Ineq. Appl. , 7(5) (2002), 673-699.
\item F.Marcellán, B.P. Osilenker, I.A.Rocha, On Fourier series of a discrete Jacobi-Sobolev Inner Product, J. Approx. Theory, 117(2002), 1-22.
\item F. Marcellán, Y.Xu, On Sobolev orthogonal polynomials, Expositiones Math., 33(2015), 308-352.
\item B.P.Osilenker, Generalized trace formula and asymptotics of the averaged Turan determinant for orthogonal polynomials, J. Approx. Theory. 141(2005), 70-94.
\item B.P.Osilenker, An Extremal Problem for Algebraic Polynomials in the Symmetric Discrete Gegenbauer-Sobolev Space, Math.Notes, 82:3(2007), 411-425; translation in Mathematical Notes, 82:3(2007), 366-379.
\item B.P.Osilenker, On linear summability methods of Fourier series in polynomials orthogonal in a discrete Sobolev space, Siberian  Math. Journal 56:2(2015), 420-435; translation in Siberian Math. Journal, 56:2(2015), 339-351.
\item B.P.Osilenker, ?n multipliers for Fourier series in polynomial orthogonal in continual-discrete Sobolev spaces, Contempo\-rary problems of mathematical and mechanics, ?oscow,Max Press, 2019, 500-503 [in Russian].
\item I.A.Rocha ,F. Marcellán , L.Salto , Relative asymptotics and Fourier series of orthogonal polynomials with a discrete Sobolev  inner product, J. Approx. Theory,121( 2003), 336-356.
\item J.M.Rodriguez, Approximation by polynomials and smooth functions in Sobolev spaces with respect to measures, J. Approx. Theory 120(2003), 185-216.
\item J.M.Rodriguez, V.Álvarez, E.Romera, D.Pestana, Generalized Weighted Sobolev Spaces   and Applications to Sobolev Orthogonal Polynomials, I,Acta Appl. Math. 80(2004), 273-308.
\item J.M.Rodriguez, V.Álvarez, E.Romera, D.Pestana, Generalized Weighted Sobolev Spaces   and Applications to Sobolev Orthogonal Polynomials, II, Approx . Theory Appl. 18:2(2002), 1-32.
\item I.I.Sharapudinov, Systems of Sobolev–orthogonal functions associated with an orthogonal systems, Izv.Ross.Akad.Nauk, Ser.Math., 82:1(2018), 225-258; translation in Izvestia Mathematics,82:1(2018), 212-244.
\item I.I.Sharapudinov. Systems of functions, Systems of Sobolev-orthogonal functions and their Applications, Uspehi Math. Nauk, 74:4(448)(2019), 87-164 [in Russian].
\item H.Triebel.Interpolation Theory,Function Spaces, Differential Operators, Mir, Moscow, 1980 [in Russian].
\end{enumerate}

\end{document}